\documentstyle[12pt,amssymb]{article}
\newcommand{\C}{\Bbb{C}}

\newcommand{\Pe}{\Bbb{P}}

\newcommand{\Z}{\Bbb{Z}}

\newcommand{{\1}}{{\bf 1}}
\newcommand{\co}{{\cal{O}}}

\newcommand{\ce}{{\cal{E}}}

\newcommand{\cf}{{\cal{F}}}
\newcommand{\cl}{{\cal{L}}}
\newcommand{\ch}{{\cal{H}}}
\newcommand{\om}{\omega}
\newcommand{\ra}{\rightarrow}
\newcommand{\lra}{\longrightarrow}
\newcommand{\ck}{{\cal{K}}}
\newcommand{\codim}{\mbox{codim}}
\newcommand{\hra}{\hookrightarrow}

\newcommand{\vp}{\varphi}

\newcommand{\pic}{\mbox{Pic}}

\newcommand{\ptp}{\Bbb{P}_1 \times \Bbb{P}_3}
\newcommand{\wed}{\widetilde{D}}
\def \D{ \displaystyle}
\begin{document}
\begin{center}
{\Large {\bf A Vector Bundle of Rank 2 on ${\bf P_1 \times P_3}$}}\\[2ex]
{\large H. Lange}
\end{center}
\vspace{3cm}
{\bf Introduction:} It is well--known that the Horrocks--Mumford bundle
$\cf$ encodes a lot of very interesting geometric information. This is
essentially the reason for the fact that  much work has been done in
order to find other rank--2 bundles on $\Pe_4$. The only non--splitting
vector bundles of rank 2 on $\Pe_4$, known upto now, are twists of
pullbacks of $\cf$ by finite coverings $f : \Pe_4 \ra \Pe_4$. So it  seems
to be a natural question to consider, instead of  $\Pe_4$, other Fano
4--folds.
 It is the aim of this
note to give an example of a rank--2 vector  bundle on $\Pe_1 \times \Pe_3$ and
to
show that it also admits very interesting  geometric properties.\\
\par
Consider $P = \Pe_1 \times \Pe_3$ and let $p$ and $q$ denote its
projections onto $\Pe_1$ and $\Pe_3$, respectively. There is a canonical
isomorphism $\Z \times \Z \ra \pic(P)$, given by $(a,b) \mapsto \co_P
(a,b)  := p^{\ast} \co_{\Pe_1} (a) \otimes q^{\ast} \co_{\Pe_3} (b)$.
First note that it is easy to construct nonsplitting vector  bundles of
rank 2 on $P$: According to the K\"unneth formula and Serre duality
$$
h^1 (P, \co_P (a,b)) = - (a+1) {{b+3}\choose {3}} > 0 $$
for all $a \leq -2$
and $b \geq 0$. Hence  for any such pair $(a,b)$ there is a rank--2 vector
bundle $\cf$ on $P$ fitting into an exact sequence
$$
0 \lra \co_P \lra \cf \lra \co_P  (-a, -b) \lra 0 .
$$
Comparing  Chern classes it is easy to see that $\cf$ is not a direct
sum of line bundles on $P$. From a geometric point of view these vector
bundles are not so interesting, for example they have no jumping lines.\\
\par
One of the reasons for the geometric relevance of the  Horrocks--Mumford
bundle seems to be that it is related, via the Serre construction, to a
nonregular surface in $\Pe_4$. This suggests that a nonregular
subcanonical locally complete intersection surface in $\Pe_1 \times
\Pe_3$ might yield an interesting vector bundle.\\
\par
In [L] the abelian surfaces in $\Pe_1 \times \Pe_3$ were classified. In
fact there is a twodimensional family of abelian surfaces  admitting an
embedding into $\Pe_1 \times \Pe_3$ and for every such surface $X$  there
are only finitely many  embeddings $X \hra \Pe_1 \times \Pe_3$.
Via the Serre construction every such embedding yields a vector bundle
$\ce$ of rank 2 with Chern classes $c_1 (\ce) = 2h_1 + 4 h_3$ and $
c_2(\ce) = 8h_1 h_3 + 6 h^2_3$, where $h_i$ denotes the pullback of the
class of a
hyperplane section on $\Pe_i$. It is easy to see that the bundle $\ce$
dannot be derived from the Horrocks--Mumford bundle on $\Pe_4$ by an
elementary transformation. The aim of this note is to derive some
geometric properties of the bundle $\ce$ to $P$.
\medskip

\noindent
Section 1 recalls  some properties of the embedding $X \hra
\Pe_1 \times \Pe_3$ from [L] and derives some additional facts needed
subsequently. In Section 2 the group of symmetries of $X$ in $\Pe_1
\times \Pe_3$ is computed.  In Section 3 it is shown that $\ce$ is a
stable vector bundle with respect to $(h_1, h_3)$. Section 4
studies  the restrictions $\ce | \{ t \} \times \Pe_3$ as a family of
vector bundles on $\Pe_3$. In Section 5 the degree of a  rational pencil
of quadrics, which can be associated to $\ce$ in a natural way, is
determined. In Section 6 the generic splitting types
and the set of jumping lines of $\ce$ are computed. Finally Section 7
contains a list of open problems.
\medskip

\noindent
I would like to thank Ch. Birkenhake and K. Hulek for some valuable
discussions. In particular Proposition 5.4 is due to the latter one.
\section{Abelian surfaces in ${\bf P_1   \times P_3}$}
 It is shown in [L]
that there is a two--parameter family of abelian surfaces
admitting an embedding into $\ptp$. Moreover, it is proven that every
abelian surface in $\ptp$ is a member of this  family. Suppose $\vp =$
$(\vp_1, \vp_3) : X \hra P := \ptp$ is such an embedding of an abelian
surface $X$ over the field of complex numbers. In this section we recall
some of its properties from [L] and derive some additional facts, which
are needed subsequently.\\

 \par
Let $h_i$ denote the pullback of the hyperplane section class of $\Pe_i$
to $P$ for $i = 1$ and $3$.\\

\par
{\bf (1)} the class $[X]$ of  $X$ in $H^4(P, \Z)$ is
$$
[X] = 8h_1 h_3 + 6h^2_3 .
$$
\par {\bf (2)} There is a commutative diagram
$$
\begin{array}{rcl}
0 \ra E \ra X&\stackrel{q}{\lra} & F \ra 0\\[1ex]
& \vp_1 \! \searrow & \downarrow \overline{\vp_1}\\[1ex]
&&\Pe_1
\end{array}
$$
where the row is an
exact sequence of abelian varieties with elliptic curves $E$ and $F$ and
$\overline{\vp}_1$ is a double covering.\\

\par
{\bf (3)} Let $L_i$ denote the line bundle on $X$ defining the morphism
$\vp_i : X \ra \Pe_i$. Then $\vp_i$ is given by the complete linear
system
$|L_i|$.\\

\par
{\bf (4)} The line bundle $L_3$ is of type $(1,4)$ on $X$ and $\vp_3 : X
\ra \Pe_3$ is a birational morphism onto a singular  octic in $\Pe_3$.\\

\par
{\bf (5)} For every $t \in \Pe_1$ consider the scheme theoretical
intersection $X_t := X \cap \{ t \} \times \Pe_3$ as a curve in $\Pe_3$,
i.e.
$$
X_t = \vp_3 (\vp^{-1}_1 (t)) .
$$
Let $t_1, \ldots , t_4$ denote the 4 ramification points of the double
covering $\overline{\vp_1} : F \ra \Pe_1$.
\bigskip

\noindent
{\bf Lemma 1.1.} (i) {\it For every $t \not= t_1, \ldots , t_4$ the
curve $X_t$ is a disjoint  union of \, $2$ smooth plane cubics $E_t$ and
$E'_t$ in $\Pe_3$ both isomorphic to $E$.\\}
(ii) {\it For $\nu = 1, \ldots, 4$ the curve $X_{t_\nu}$ is a double curve in
$\Pe_3$, not lying in a plane, with support a plane cubic $E_{t_\nu}$
isomorphic to } $E$.
\bigskip

\noindent
{\bf Proof. } (i) This follows from (2) and (3), since $\vp :X \hra P$
embeds the curve $X_t$  into $\{ t \} \times
\Pe_3$.\\
(ii) The assertion about the support follows again from (2) and (3).
That $X_{t_{\nu}}$ is a double curve follows from (2) and the fact that
$t_{\nu}$ is a ramification point of the double covering $\overline{\vp_1}$.
It is only to show that the linear system $|L_3| \big| X_{t_{\nu}}$is of
dimension 3.\\
Since dim$|L_3|=3$, we have obviously dim$(|L_3|\big| X_{t_v}) \leq 3$.
For the converse inequality consider $X$ as a family of curves $(X_t)_{t
\in \Pe_1}$ over $\Pe_1$. Since dim$(|L_3| \big|X_t) =3$ for a general $t
\in \Pe_1$, the assertion follows by semicontinuity. \hfill $\Box$\\

\par
{\bf (6)} For every point $x \in F$ let $P_{2,x}$  denote the plane in
$\Pe_3$ spanned by the plane cubic $\vp_3(q^{-1}(x)) \,\, ( = E_t $ or $
E'_t$ if $t =  \overline{\vp}_3 (x)).$
\bigskip

\noindent
{\bf Lemma 1.2.} $P_{2,x} \not= P_{2,y}$ {\it for all} $x,y \in F, x \not=
y$.
\bigskip

\noindent
{\bf Proof.} Suppose $P_{2,x} = P_{2,y}$. Let $D \in |L_3|$ denote the
divisor corresponding to the plane $P_{2,x} = P_{2,y} $ in $\Pe_3$. then
$$
D =q^{-1} (x) + q^{-1} (y) + \widetilde{D}
$$
with elliptic curves $q^{-1}(x)$ and $q^{-1} (y)$ and a curve
$\widetilde{D}$.
It follows
$$
8 = D^2 = D \cdot q^{-1} (x) + D \cdot q^{-1} (y) + D \cdot
\widetilde{D} = 6 + D \cdot \widetilde{D}
$$
and hence
$$
2 = D \cdot \wed = q^{-1} (x) \wed + q^{-1} (y) \wed + \wed^2 .
$$
Since all summands on the right hand side are nonnegative, one of the
summands must be 0. But then $\wed$ is a disjoint union of translates  of
$E$ and all summands are $0$, a contradiction. \hfill $\Box$\\

\par
{\bf
(7)} Let $x' \in F$ denote the conjugate of $x \in F$ with respect
to the double covering $\overline{\vp_1}$. Then, denoting $t =
\overline{\vp_1} (x)$,
$$
Q_t := P_{2,x} + P_{2,x'}
$$
is a quadric in $\Pe_3$, of rank 2 for all $t \not= t_1, \ldots, t_4$ and
of rank 1 for $t = t_\nu$. Consider the pencil $\{ Q_t | t \in \Pe_1 \}$
of these quadrics.
\bigskip

\noindent
{\bf Lemma 1.3.} $\{Q_t | t \in \Pe_1 \}$ {\it is a rational pencil of
quadrics of rank $\leq 2$ and degree $ d\geq 2$ in $\Pe_3$, i.e. of the
form $\{ \lambda^d Q_0 + \lambda^{d-1} Q_1 +\ldots + \mu^d Q_d | (\lambda,
\mu) \in \Pe_1 \}$ with quadrics $Q_0, \ldots, Q_d $ in }$\Pe_3$.
\bigskip

\noindent
{\bf Proof.} Suppose that $\{ Q_t | t \in \Pe_1 \}$ is a linear pencil
and let $Q$ and $Q'$ denote 2 different quadrics of rank 2 in it. We may
choose the coordinates of $\Pe_3$ such that $Q$ is given by the matrix
diag$(1,1,0,0)$. Then $Q'$ is given by a matrix $
\left(
\begin{array}{cccc}
a_0 & a_1 &0&0\\
a_1&a_2 &0&0\\
0&0&0&0\\
0&0&0&0
\end{array}
\right)$,
since all quadrics $Q_t$ are of rank $\leq 2$. Then the linear system is
given by the matrices\\
$$
\left\{ \left( \begin{array}{cccc}
\lambda a_0 + \mu &  \lambda a_1 &0 & 0\\
\lambda a_1 & \lambda a_2 + \mu &0 & 0\\
0&0&0&0\\
0&0&0&0
\end{array}
\right) i (\lambda, \mu) \in \Pe_1 \right\}
.$$
Since det${{\lambda a_0 + \mu \quad \lambda a_1}\choose{\lambda a_1  \quad
\lambda a_2 + \mu}} =0$ has only 2 solutions, the pencil admits only 2
quadrics of rank 1, a contradiction. \hfill $\Box$
\bigskip

\noindent
For a quadric $Q_t$ of rank 2 in $\Pe_3$ let $\ell_t$ denote the
singular line of $Q_t$.
\bigskip

\noindent
{\bf Lemma 1.4.} {\it The lines $\ell_t, t \in \Pe_1$ form a
one--dimensional family of lines in } $\Pe_3$.
\bigskip

\noindent
{\bf Proof.} If this is not the case, then all quadrics of the pencil
contain the same singular line $\ell$ and the pencil $\{Q_t | t \in
\Pe_1 \}$ is given by a curve $q$ of degree $d \geq 2$ in the plane $\pi$
parametrizing the quadrics in $\Pe_3$ with singular line $\ell$. On the
other hand the quadrics in $\pi$ containing a fixed plane in $\Pe_3$ form
a line in $\pi$. This line intersects $q$ in $d \geq 2$ points, which
contradicts Lemma 1.2. \hfill $\Box$
\medskip

\noindent
We will see later (see Proposition 5.3) that the pencil $\{ Q_t | t \in \Pe_1
\}$ is of degree 4.
\section{Symmetries of $X$ in ${\bf P_1 \times P_3}$. } Let $\vp = (\vp_1,
\vp_3)
:X \ra  \Pe_1 \times \Pe_3$ be an embedding of an abelian surface $X$
defined by line bundles $L_1$ and $L_3$ on $X$ as in Section 1. In this section
we
want  to determine the group of all translations of $X$ which extend to
elements of $P G L_1(\C) \times P G L_3 (\C)$.\\

\par
Recall that for any line bundle $L$ on $X$ the subgroup of all $x \in
X$ with $t^{\ast}_x L \simeq L$ is denoted by $K(L)$. Here $t_x : X \ra
X$ denotes the translation of the abelian variety $X$ by $x$. If $L$ is base
point free with
corresponding morphism $\vp_L : X \ra \Pe_m$, then $K(L)$ is the
largest
subgroup of translations of $X$ which are induced by automorphisms of
$\Pe_m$. If $L'$ is a second base point free line bundle on $X$ with
corresponding morphism $\vp_{L'} : X \ra \Pe_n$, let $K(L, L')$ denote
the subgroup of all translations of $X$ which are  induced by elements of
$P GL_m(\C) \times P G L_n (\C)$. Clearly such an element induces a
translation $t_x$ on $X$  if and only if $x \in K(L) \cap K(L')$.
\medskip

\noindent
Let  $X_2$ denotes the subgroup of 2--division points of $X$. Then we have
\bigskip

\noindent
{\bf Proposition 2.1.} $\quad K(L_1, L_3) = K(L_3) \cap X_2 .
$
\bigskip

\noindent
{\bf Proof.}
By construction of $\vp_1 : X \ra \Pe_1$ there is a line bundle $\ell_1$
on the elliptic curve $F$ such that $L_1 = q^{\ast} \ell_1$ with $q : X
\ra F$ as in (2). We have $K(\ell_1) = F_2$, since $\ell_1$ is of degree
2 on $F$. Hence $K(L_1) = q^{\ast} K(\ell_1)$ contains $X_2$, and thus
$$
K(L_3) \cap X_2 \subseteq  K(L_3) \cap K(L_1) = K(L_1, L_3).
$$
$L_3$ being of type $(1,4)$ implies $K(L_3) \simeq \Z / 4 \Z \times \Z /
4 \Z$. Hence it suffices to show that no element of order 4 in $K(L_3)$
is contained in $K(L_1)$. Suppose $x \in K(L_1) \cap K(L_3)$ is of order
4. Then $t^{\ast}_x (L_1 \otimes L_3) \simeq t^{\ast}_x L_1 \otimes t^{\ast}_x
L_3 \simeq L_1
\otimes L_3$, i.e. $x \in K(L_1 \otimes L_3)$. But
$$
(L_1 \otimes L_3)^2 = 2(L_1 \cdot L_3) + (L^2_3)= 20
$$
since $(L_1 \cdot L_3) =6$ (see [L], Lemma 1.1). Hence $L_1 \otimes L_3$ is of
type $(1,10)$ and $K(L_1 \otimes L_3) \simeq \Z/ 10 \Z \otimes \Z / 10
\Z$. Hence $K(L_1 \otimes L_3)$ does not contain an element of order 4.
\hfill$\Box$
\medskip

\noindent
Since $L_3$ is of type $(1,4)$, Proposition 2.1 implies
$$
K(L_1, L_3) \simeq \Z / 2 \Z \times \Z / 2 \Z
$$
and one may choose coodinates of $\Pe_1$ and $\Pe_3$ in such a way (see
[CAV], p. 169) that generators $\sigma$ and $\tau$ of $K(L_1, L_3)$ are
given by
$$
\begin{array}{ccl}
\sigma &  = &  \left( \left( \begin{array}{cc}
0&1\\
1&0
\end{array}
\right)
, \left(
\begin{array}{cccc}
0& 0& 1& 0\\
0& 0& 0&1\\
1 & 0& 0& 0\\
0 & 1 & 0& 0
\end{array}
\right)
\right) \in G L_2 (\C) \times G L_4 (\C)\\[2ex]
\tau & = & \left( \left( \begin{array}{rr}
1&0\\
0 & -1
\end{array} \right), \left(
\begin{array}{ccrr}
1& 0 & 0&0\\
0 & -1 & 0 & 0\\
0 & 0 & 1& 0\\
0 &0& 0& -1
\end{array}
\right) \right) \in G L_2 (\C) \times G L_4 ( \C).
\end{array}
$$
The {\it Heisenberg group of the pair} $(L_1, L_3)$ is by definition the
subgroup $H(L_1, L_3)$ of $G L_2(\C) \times GL_4(\C)$ generated by
$\sigma$ and $\tau$. The Heisenberg group $H(L_1, L_3)$ acts in the
usual way on the pair of line bundles $(L_1, L_3)$.
\bigskip

\noindent
{\bf Proposition 2.2.} {\it The Heisenberg group $H(L_1, L_3)$ is the
dihedral
group $D_8$ with $8$ elements. It fits into an exact sequence}
$$
0 \ra \mu_2 \ra H(L_1, L_3) \ra K(L_1, L_3) \ra 0
$$
{\it where $\mu_2$ denotes the group of order}   2 {\it generated by the
commutator } $\sigma \tau \sigma^{-1} \tau^{-1}$.
\bigskip

\noindent
{\bf Proof:} $H_1 (L_1, L_2)$ is generated bt $\tau$ and the element
$\kappa   := \sigma \tau$ of order 4 with relation $\tau \kappa  \tau =
\kappa^{-1}$.
\hfill $\Box$
\section{The vector bundle $\ce$ on ${\bf P_1 \times P_3}$}
Let $\vp = (\vp_1, \vp_3) : X \ra P = \Pe_1 \times \Pe_3$ be an
embedding as in Section 1. Denote $
\cl = \co_P(2,4)$. Then $\cl =  \om^{-1}_P$ and there is an isomorphism
$$
\xi : \cl \otimes \om_P \otimes \co_X \ra \om_X = \co_X .
$$
According to the Serre--construction (which applies in this case, see
[H], Remark 1.1.1) the pair $(X, \xi)$ determines uniquely a triplet
$(\ce, s, \psi)$ with a rank--2 vector bundle $\ce$ on $P$, a section
$ s \in H^0 (P,
\ce)$ such that $(s) = X$ and an isomorphism $\psi : \bigwedge^2 \ce \ra
\cl$.\\

\par
The Koszul complex of the section $s$ is
$$
0 \ra \bigwedge^2 \ce^{\ast} \ra \ce^{\ast} \stackrel{s}{\ra} I_X \ra 0
$$
where $I_X$ denotes the ideal sheaf of $X$ in $P$. Tensoring with
$\bigwedge^2 \ce = \cl = \co_P (2,4)$ and using $\ce =\ce^{\ast}(\bigwedge^2
\ce)$ we obtain the exact sequence
$$
 0 \ra \co_P \stackrel{s}{\ra} \ce \ra I_X (2,4) \ra 0 . \eqno(1)
 $$
Then we have, using Section 1 (1)
$$
c_1 (\ce) = 2 h_1 + 4 h_3 \eqno(2)
$$
and
$$
c_2 (\ce) = 8 h_1 h_3 + 6h^2_3 .\eqno(3)
$$

\par
Consider the ample line bundle $\ch = \co_P(1,1)$ on $P$. A vector bundle
$\cf$ of rank 2 on $P$ is called ({\it semi--) stable with respect to}
$\ch$ if for every invertible subsheaf $\cl $ of $\cf$
$$
\cl \cdot \ch^3 \stackrel{<}{(=)} \frac{1}{2} \det \cf \cdot \ch^3 .
$$
\bigskip

\noindent
{\bf Proposition 3.1.} $\ce$ {\it is stable with respect to} $\ch$.\\

\par
For the proof we need the following Lemma
\bigskip

\noindent
{\bf Lemma 3.2.} (i) {\it Suppose  $h^0 (I_X (a,b)) \not= 0$ for some
integers $a$ and $b$. Then $a \geq 0$ and} $b \geq 2$.
\begin{itemize}
\item[(ii)] $h^0 (I_X (0,b)) \not= 0$ {\it if and only if} $b \geq 8$.
\item[(iii)] $h^0(I_X(a,2)) = 0$ {\it for }$a \leq 2.$
\end{itemize}
\bigskip

\noindent
{\bf Proof.} (i)  Since $h^0( \co_P (a,b)) =0$ if $a < 0$, it suffices to
show that $b \geq 2$. For a general $t \in \Pe_1$ the intersection
$X_t = X \cap \{ t \} \times \Pe_3$ is a disjoint untion of 2 plane
cubics in $\Pe_2$ (see Section 1 (5)).
In particular $h^0(\Pe_3, I_{X_t} (1)) =0$. So the exact sequence
$$
 0 \ra I_X (a,1) \ra I_{X_t} (1) \ra I_{X_t} (1) / I_X (a,1) \ra 0
 $$
 implies $h^0 (I_X (a,1)) =0$ for all $a$.\\
  (ii ): $h^0 (I_X (0,b)) \not= 0$ means that $X \subset \Pe_1
 \times V_b$ where $V_b$ is a surface of degree $b$ in $\Pe_3$. Now that
 the composed map $\vp_3 : X \hra \Pe_1 \times \Pe_3 \ra \Pe_3$ is
 birational onto a  (singular) octic in $\Pe_3$ (see  Section 1(4)). This
implies
 the assertion.\\
  (iii): It suffices to show $h^0(I_X (2,2))= 0$. But $h^0(I_X
 (2,2)) \not= 0$ would mean (using the results of Section 4)  that the pencil
$\{ Q_t | t \in \Pe_1 \}$ of
 Section 1 (7) is quadratic, contradicting Proposition 5.3  below. \hfill
$\Box$
\bigskip

\noindent
{\bf Proof of Proposition 3.1.} Suppose $\ce$ is not stable with respect
to $\ch$ and let $\co_P(a,b)$ be a subinvertible sheaf of $\ce$
violating the stability of $\ce$, i.e.
$$
a + 3b = \co_P (a,b) \cdot \ch^3 \geq \frac{1}{2} \det \ce \cdot \ch^3 =
7 .
$$
Then the composed map
$$
\co_P(a,b) \ra \ce \ra I_X (2,4)
$$
is nonzero,
implying $$
h^0(I_X(2-a, 4-b)) \not= 0.
$$
Lemma 3.2 implies that either $a
\leq 1$ and $b \leq 2$ and $a \leq -1$ if $b=2$ or $a=2$ and $b \leq
-4$. But in the first case $7 \leq a+3b \leq 5$ and in the second case $7
\leq a + 3b \leq - 10$, a contradiction. \hfill$\Box$
\bigskip

\noindent
{\bf Remark 3.3.} The line bundle $\co_P (m,n)$ is ample if and only if
$m > 0$ and $n > 0$. Defining the (semi--) stability of a rank--2 vector
bundle on $P$ with respect to $\co_P(m,n)$ in the same way as for $\ch$.
Then an analogous proof yields that $\ce$ is stable (respectively
semistable) with respect to $\co_P (m,n)$ if and only if $n < 18 m $
(respectively $n \leq 18 m)$.
\medskip

\noindent
 In order to lift the action of the Heisenberg group $H(L_1, L_3)$ of the
 pair $(L_1, L_3)$ to the vector bundle $\ce$, we need the cohomology of
 $\ce (-2, -4)$.
 \bigskip

 \noindent
 {\bf Lemma 3.4} $$
 h^i(\ce(-2, -4)) = \left\{ \begin{array}{lcr}
 0 && i \not= 2\\
 &\mbox{\it for} &\\
 2 && i=2
 \end{array}
 \right.
 .$$
 \medskip

 \noindent
 {\bf Proof:} Serre duality says $h^{4-i} (\ce(-2, -4)) = h^i (\ce(-2,
 -4))$. Hence it suffices to compute $h^i(\ce ( -2, -4))$ for $i \leq
 2$. The exact sequence (1) yields
 $$
 h^i( \ce (-2, -4)) = h^i(I_X)
 $$
 for $i = 0,1,2$. Now the exact sequence $0 \ra I_X   \ra \co_P \ra \co_X \ra
 0$ implies the assertion. \hfill $\Box$
 \bigskip

 \noindent
 {\bf Proposition 3.5:} {\it The vector bundle $\ce$ admits an action of
 the Heisenberg group $H(L_1, L_3)$ uniquely determined upto a constant
 by the embedding $X \hookrightarrow P$.}
 \bigskip

 \noindent
 {\bf Proof:} An element $\mu \in H(L_1, L_3)$ may be considered as an
 automorphism of $P$. Moreover $\mu$ defines an isomorphism $\vp_{\mu} :
 I_X (2,4) \widetilde{\ra} \mu^{\ast} I_X(2,4)$. Consider the diagram
 $$
 \begin{array}{llcccccr}
 0 \ra & \co_P & \stackrel{s}{\ra} &\ce& \ra & I_X(2,4)& \ra 0\\
 & ||& &&& \downarrow \vp_{\mu} &\\
 0 \ra & \co_P & \ra & \mu^{\ast} \ce& \ra & \mu^{\ast} I_X(2,4) & \ra 0
 \end{array}
 $$
 It induces an exact sequence
 $$
 0 \ra \, \mbox{Hom} (\ce, \mu^{\ast} \ce) \ra \, \mbox{Hom} (\ce,
 \mu^{\ast} I_X(2,4)) \ra \, \mbox{Ext}^1(\ce, \co_P)
 $$
 But Ext$^1(\ce, \co_P) = H^1(\ce^{\ast}) = H^1(\ce(-2, -4))=0$
 according to Lemma 3.4. This implies the assertion, the uniqueness
 coming from the fact that the embedding $X \hookrightarrow P$
 determines the section $s$ upto a constant. \hfill $\Box$
 \section{Restriction to ${\bf \{ t \} \times P_3}$.}
For any point $t \in \Pe_1$ consider the restriction
$$
\ce_t := \ce | \{ t \} \times \Pe_3
$$
as a rank--2 vector bundle on $\Pe_3$. \,\, $X_t$ is a local complete
intersection curve in $\Pe_3$ and the triplet $(\ce_t, s | \{ t \}
\times \Pe_3, \psi | \{ t \} \times \Pe_3 )$ corresponds to the pair
$(X_t, \xi | \{ t \} \times \Pe_3)$ via the Serre--construction. Hence
we get the following exact sequence directly by restricting
(1) of Section 2.
$$
0 \ra \co_{\Pe_3} \stackrel{s_t}{\ra} \ce_t \ra I_{X_t} (4) \ra 0. \eqno(1)
$$
This implies
$$
\begin{array}{l}
c_1 (\ce_t) =4\\
c_2 (\ce_t) = 6
\end{array}
$$
\medskip

\noindent
{\bf Proposition 4.1} {\it For any $t \in \Pe_1 $ the vector bundle
$\ce_t$ is semistable but not stable.}
\vspace{1cm}
\newline
{\bf Proof.} The quadric $Q_t$ of Lemma 1.3 contains the curve $X_t$ and
is in fact the only quadric in $\Pe_3$ containing $X_t$. Hence $h^0
(I_{X_t} (2)) =1$ and the exact sequence $0 \ra \co_{\Pe_3} (-2) \ra E_t
(-2) \ra I_{X_t} (2) \ra  0$ implies $h^0 (E_t(-2)) =1.$ This means that
$E_t$ is not stable. On the other hand $h^0(E_t(-3)) =  0$ means   that
$E_t$ is semistable.\hfill $\Box$\\

\par
Let $\sigma_t \in H^0 (E_t(-2))$ denote the non--zero section, uniquely
determined upto a constant, and $Y_t = (\sigma_t)$ the corresponding
zero variety in $\Pe_3$. Then we have the exact sequence
$$
0 \ra \co_{\Pe_3} \stackrel{\sigma_t}{\ra} \ce_t (-2) \ra I_{Y_t} \ra 0
. \eqno(2)
$$
\smallskip

\noindent
{\bf Lemma 4.2.} $Y_t$ {\it is a curve of degree $2$ in $\Pe_3$ with
$p_a (Y) = -3$ and} $\omega_{Y_t} = \co_{Y_t} (-4)$.
\vspace{1cm}
\newline
{\bf Proof.} According to [H], Proposition 2.1 we have
$\deg Y = c_2 (E(-2))=2$ and $2p_a(Y)-2= c_2(E(-2)) \cdot (c_1(E(-2))-4) =
- 8$,
i.e. $p_a(Y) = -3$. The assertion for $\omega_{Y_t}$ follows from the
adjunction formula.\hfill $\Box$
\vspace{0.51cm}

\par
Let $Z_t$ denote the reduced curve underlying $Y_t$.\\
\bigskip

\noindent
{\bf Lemma 4.3} (i): \,\, $Z_t$ {\it is a line in $\Pe_3$ and $Y_t$ is a
multiplicity } 2 {\it structure on} $Z_t$.
\\
(ii): {\it Choose the coordinates $x_i , \quad i = 0, \ldots , 3 $ of $\Pe_3$
in
such a way that $Z_t$ is the line $x_0 = x_1 =0.$ then $Y_t$ is given by
the homogeneous ideal}
$$
(x_0^2, x_0 x_1, x^2_1, f \cdot x_0 + g \cdot x_1)
$$
{\it where $f$ and $g$ are forms of degree $3$ in $x_2$ and $x_3$ without
common zero}.\\
\bigskip

\noindent
{\bf Proof.}
 (i): Since $\deg Y_t =2$, it suffices to show that $Y_t$ is irreducible.
 But otherwise $Y_t$ would be reduced and hence $p_a(Y) \geq -1$.\\
 (ii):  According to Ferrand's theorem (see [H], Theorem 1.5 with $m=4$)
 there is an exact sequence
 $$
 0 \ra I_{Y_t} \ra I_{Z_t} \ra \co_{Z_t} (2) \ra 0.
 $$
Consider the conormal bundle sequence $0 \ra N^{\ast}_{Y_t|\Pe_3} |_{ Z_t}
\ra N^{\ast}_{Z_t|\Pe_3} \ra N^{\ast}_{Z_t|Y_t} \ra 0$. According to
(i)
\quad $N^{\ast}_{Z_t|\Pe_3} = \co_{\Pe_1} (1)^{\oplus 2}$.
 Hence we get an exact sequence
 $$
 0 \ra I_{Y_t} | \Pe_1 \ra \co_{\Pe_1} (1)^{\oplus 2}
 \stackrel{u}{\ra} \co_{\Pe_1} (2) \ra 0 .
 $$
 The map $u$ is given by 2 forms $f$ and $g$ of degree 3 on
 $\Pe_1$ without common zeros. This implies that $I_{Y_t} | \Pe_1$ is
 generated by $f x_0 + g x_1$.\hfill $\Box$\\

\par
Suppose $t \not= t_{\nu}$ for $\nu = 1, \ldots , 4$. Then $X_t$ consists
of 2 disjoint plane cubics $E_t$ and $E'_t$. Denote by $P_t$ and $P'_t$
the planes spanned by $E_t$ and $E'_t$. Then we have\\
\bigskip

\noindent
{\bf Proposition 4.4.} {\it The line $Z_t$ is the line of intersection
$P_t \cap P'_t$ and the cubic forms $f$ and $g$ are given by $E_t \cap
Z_t$ and } $E'_t \cap Z_t$.\\
\bigskip

\noindent
{\bf Proof.} Restricting the exact sequence (1) to the plane $P_t$ we
obtain an exact sequence
$$
0 \ra \co_{P_t} (3) \ra \ce_t | P_t \ra I_{E'_t \cap P_t} (1) \ra 0.
\eqno(3)
$$
since $P_t$ contains the plane cubic $E_t$. Computing the generic
splitting type of $\ce_t | P_t$ using (2), we see that $P_t$ contains
the line $Z_t$ (but not the double curve $Y_t)$. Similarly $P'_t$
contains $Z_t$, i.e. $P_t \cap P'_t = Z_t$. Restricting (2) to $P_t$ and
comparing it with (3) gives the assertion for the forms $f$ and $g$.
\hfill $\Box$
\bigskip

\noindent
\section{The Double Structure and the Pencil of Quadrics}
The quadrics of the rational pencil $\{Q_t  | t \in \Pe_1 \}$ of Section
1 fill up a 3--fold
$Q \in \Pe_1 \times \Pe_3$ containing $X$. According to Lemma 1.3
\, \, $\{
Q_t| t \in \Pe_1 \}$ is a rational pencil of degree $d \geq 2$ in
$\Pe_3$. Hence $Q$ is a  hypersurface of bidegree $(d,2)$ in $\Pe_1
\times \Pe_3$. Obviously $Q$ is the  only hypersurface of bidegree
$(d,2)$ in $\Pe_1 \times \Pe_3$ containing $X$. This means
$$
h^0 (I_X (d, 2)) =1.
$$
 From the exact sequence
$
0 \ra \co_{\Pe_1 \times \Pe_3} (d-2, -2) \ra \ce (d-2, -2) \ra I_X(d,2)
\ra 0
$ we deduce
$$
h^0(\ce (d-2, -2)) =1.
$$
Hence there is an exact sequence,
unique upto a multiplicative constant
$$
0 \ra \co_{\Pe_1 \times \Pe_3} \stackrel{\sigma}{\ra} \ce(d-2, -2) \ra
I_Y(2d-2, 0) \ra 0 \eqno(1)
$$
with a surface $Y \subseteq \Pe_1 \times \Pe_3$ of class
$$
[Y] = 4h_1 h_3 + 2h^2_3.
$$
 For every point  $t \in \Pe_1$ the exact sequence (1) restricts to the exact
 sequence (2) of Section 4. Hence for every $t \in \Pe_1$ the surface $Y
 \subseteq \Pe_1 \times \Pe_3$ restricts to the curve $Y_t \subseteq \{
 t \} \times \Pe_3$. According to Lemma 4.3 \, \, $Y_t$ is a double structure
on a line
 $Z_t$ in $\{ t \} \times \Pe_3$ for all $t$. This implies\\
 \bigskip

 \noindent
 {\bf Proposition 5.1.} {\it The surface $Y \subseteq \Pe_1 \times
 \Pe_3$ is a double structure on a $\Pe_1$--bundle $Z \subseteq \Pe_1
 \times \Pe_3$ over $\Pe_1$. The bundle map $Z \ra \Pe_1$ coincides with
 the projection onto the first factor of $\Pe_1 \times \Pe_3$.
 The class of $Z$ in $H^4(\Pe_1 \times \Pe_3,\Z)$ is $[Z]
 = 2h_1 h_3 + h^2_3$.}\\
 \bigskip

 \noindent
 In other words $Z$ is a certain Hirzebruch surface $
 \Sigma$ embedded into $\Pe_1 \times \Pe_3$ with class $2h_1 h_3 +
 h^2_3$ such that the bundle map  $\Sigma \ra \Pe_1$ coincides with the
 projection onto the first factor. The following lemma classifies these
 embeddings.\\

\par
 Let $\Sigma_e$ denote the Hirzebruch surface with invariant $e \geq
 0$.
 There is a section $C_0 $ of $p : \Sigma_e \ra \Pe_1$ such that
 $$
 \pic (\Sigma_e) \simeq \Z [C_0] + \Z[f]
 $$
 where $f$ denotes a fibre of $p$, with intersection numbers
 $$
 C^2_0 = -e \quad , \quad C_0 \cdot f =1 \quad , \quad f^2 =0.
 $$
 Suppose $i = (\vp_1, \vp_3): \Sigma_e \hra \Pe_1 \times \Pe_3$ is an
 embedding such that
 $$
 [\Sigma_e ] = 2h_1 h_3 + h^2_3 \,\, \mbox{in} \,\,  H^4(\Pe_1
 \times \Pe_3, \Z) \leqno(i)
 $$
 $$
  \vp_1 :
\Sigma_e \ra \Pe_1 \,\, \mbox{is given by the linear system} \quad  |f|.
\leqno(ii)
$$
Then we have
\bigskip

\noindent
{\bf Lemma 5.2.} {\it Either $e =0$ and $\vp_3$ is given by the linear
system $|C_0 +f|$ or $e = 2$  and $\vp_3$  is given by the linear
system $|C_0 + 2f|$. Conversely, in these cases $i = (\vp_1 , \vp_3)$ is
an embedding.}
\bigskip

\noindent
{\bf Proof.} As usual we identify classes in $H^i( \Sigma_e, \Z)$ with
their images in $H^i (\Pe_1 \times \Pe_3, \Z)$. Condition (ii) means
$[ \Sigma_e] \cdot h_1 = f$ which implies
$$
f \cdot h_1 =0 \quad \mbox{and} \quad f \cdot h_3 =1. \eqno(2)
$$
Let $\vp_3 : \Sigma_e  \ra \Pe_3$ be given by a sublinear system of $|
\alpha C_0 + \beta f|$.
Necessarily $\alpha \geq 0$ and $\beta \geq 0$.
$$
2h_1 h^2_3 + h^3_3 = [ \
\Sigma_3] \cdot h_3 = \alpha C_0 + \beta f . \eqno(3)
$$
Using (2) this implies
$$
\alpha C_0 \cdot h_1 =1
$$
and
$$
\alpha C_0 \cdot  h_3 = 2 - \beta.
$$
This implies
$$
C_0 \cdot h_1 = 1 \quad \mbox{and} \quad \alpha =1 .
$$
Hence
$$
\beta \leq 2 .
$$
On the other hand
$$
2 = [
\Sigma_e] \cdot h^2_3 = (h_3 |
\Sigma_e)^2 = (C_0 + \beta f)^2 = - e + 2 \beta
$$
i.e.
$$
\beta = \frac{e}{2} +1 .
$$
We remain with the 2 cases $e = 0, \, \beta =1$ and $e = 2, \, \beta = 2$.
\\
Conversely assume that the surface and the linear systems are of one of
the 2 cases.
In
the first case $\vp_3$ embeds $\Sigma_0 = \Pe_1 \times \Pe_1$ as a
smooth quadric in $\Pe_3$. In the second case $\vp_3$ embeds $\Sigma_2 -
C_0$ and contracts $C_0$ to a quadric cone in $\Pe_3$. But then $\vp_1$
is injective on $C_0$ and the differential of $(\vp_1, \vp_3)$ has maximal
rank.
 So in both cases $i = (\vp_1, \vp_3)$ is an embedding. \hfill
$\Box$
\\

\par
The surface $Y$ is a double structure on $Z$ with $Z$ as in Lemma 6.2.
Let $\cl$ denote the ideal of $Z$ in $\co_Y$. So we have the exact
sequence
$$
0 \ra \cl \ra \co_Y \ra \co_Z \ra 0 \eqno(4)
$$
Note that $\cl$ is a line bundle on $X$, since $\cl/\cl^2 =0$.
Using the Riemann--Roch Theorem for the vector bundle $\ce (a, b),
$ saying
$$
\chi( \ce (a,b)) = - 6 + 12 a + \frac{34}{3} b + 6 b^2 + \frac{41}{3} a
b + \frac{2}{3} b^3 + 4 ab^2 + \frac{1}{3} ab^3 \eqno(5)
$$
for all $a,b \in \Z$ one can determine $\cl$ and $d$:
\bigskip

\noindent
{\bf Proposition 5.3.} {\it Suppose $Y \subseteq \Pe_1 \times \Pe_3$ is
a double structure on $Z$ ($\simeq
\Sigma_e, e = 0$ or $2$) defined by the line bundle $\cl$ on $Z$. Then}
\begin{itemize}
\item[(a)] $
\cl = \left\{ \begin{array}{lcc}
\co_Z(2 C_0 - 2f) && e = 0\\
& \mbox{\it if} &\\
\co_Z(2 C_0) && e=2
\end{array}
\right.
$
\item[(b)] {\it The pencil of quadrics $\{Q_t | t \in \Pe_1 \}$ in $\Pe_3$ is
of degree} $4$.
\end{itemize}
\par
\bigskip

\noindent
{\bf Proof.} Let $\cl = \co_Z (x C_0 + yf)$ for some $x, y \in \Z$. For
$\alpha, \beta \in \Z$ we have
$$
\co_Z \otimes \co_{\Pe_1 \times \Pe_3} (\alpha, \beta) = \left \{
\begin{array}{lcc}
\co_Z (\beta C_0 +(\alpha + \beta)f) && e=0\\
& \mbox{if}&\\
\co_Z (\beta C_0 +(\alpha + 2 \beta )f) && e=2
\end{array}
\right.
$$
and from (1) and (4) we get the exact sequences
$$
0 \ra \co_{\Pe_1 \times \Pe_3} (a -  d+2, b+2) \ra \ce(a,b) \ra I_Y (a+d,
b+2) \ra 0 \eqno(6)
$$
$$
0 \ra I_Y (a+d, b+2) \ra \co_{\Pe_1 \times \Pe_3} (a+d, b+2) \ra \co_Y
(a+d, b+2) \ra  0\eqno(7)
$$
and
$$
\scriptstyle{ 0 \ra \co_Z((x+b+2) C_0 +(y+a+b+d+2)f) \ra \co_Y (a+d, b+2) \ra
\co_Z((b+2)C_0 +(a+b+d+2)f) \ra 0} \eqno(8)
$$
if $e = 0$ and
$$
\scriptstyle{ 0 \ra \co_Z ((x+b+2)C_0 +(y+a+2b+d+4)f)\ra \co_Y (a+d, b+2) \ra
\co_Z((b+2)
C_0+ (a+2b+d+4)f)  \ra 0 }\eqno(9)
$$
if $e=2$. Using Riemann--Roch for line bundles on $\Pe_1 \times \Pe_3$
and $Z$, equations (6) -- (9) yield
$$
\scriptstyle{ \chi(\ce(a,b))=\chi(\co_{\Pe_1 \times \Pe_3} (a-d+2, b+2)) +
\chi(\co_{\Pe_1 \times \Pe_3} (a+d, b+2)) - \chi(\co_Y(a+d, b+2))}
$$
$$
\scriptstyle{ =(a-d+3){{b+5}\choose{3}} +(a+d+1) {{b+5}\choose{3}} -}
$$
$$
\scriptstyle{
\left\{
\begin{array}{l}
\scriptstyle{ - \chi(\co_Z((x+b+2)C_0
+(y+a+b+d+2)f))-\chi(\co_Z((b+2) C_0+(a+b+d+2)f))
\quad  \mbox{if} \quad \, e=0}\\
\scriptstyle{ - \chi(\co_Z((x+b+2)C_0+(y+a+2b+d+4)f)) - \chi(\co_Z((b+2) C_0 +
(a+2b
+d+4)f))
\quad  \mbox{if} \,\quad  {e=2} }
\end{array}
\right.
}
$$
$$
\scriptstyle{
= (2a+4){{b+5}\choose{3}} + \left\{
\begin{array}{l}
\scriptstyle{ -(x+b+3)(y+a+b+d+3)-(b+3)(a+b+d+3) \,\quad  \mbox{if}\,\quad
e=0}\\
\scriptstyle{ -(x+b+3)(y-x+a+b+d+3)-(b+3)(a+b+d+3) \,\quad  \mbox{if} \,\quad
e=2.}
\end{array}
\right.
}
$$
Comparing this with (5) we obtain
$$
(x-2)a + (x+y+2d-8)b +(y+d+3)(x+3)+3d -37 =0 \quad \,\, \mbox{if} \,\, e=0
$$
and
$$
(x-2)a +(y+2d-8)b+(y-x+d +3)(x+3)+3d-37 =0 \quad \,\, \mbox{if} \,\, e=2
$$
as polynomials in $a$ and $b$.
\bigskip

\noindent
Hence if $e=0 : x =2, y = 6-2d$ and $2d -8 =0$ and if $e=2 : x=2, y=8
-2d$ and $2d -8 =0$. \hfill $\Box$\\

\par
K. Hulek gave a different proof of Proposition 5.3 using the normal
bundle sequence of $Z \subseteq Y \subseteq P$. The following
proposition saying that the case $Z = \Sigma_2$ does not occur in our
situation is also due to him.
\vspace{1cm}
\newline
{\bf Proposition 5.4.} (Hulek): {\it Suppose $Z ( = \Sigma_e)$ embedded
in $P = \Pe_1 \times \Pe_3$ as above admits a double structure $Y$ in
$P$, then} $ e = 0$.
\bigskip

\noindent
{\bf Proof:} Suppose $Z = \Sigma_2 \subset P$ admits a double structure
$Y \subseteq P$.
According to Lemma 5.2 the morphism $\vp_3 | Z : Z \ra \Pe_3$ is the
contraction of the curve $C_0$ with image a quadric in $\Pe_3$. Hence
$C_0$ considered as a curve in $P$ is the preimage of a point in $\Pe_3$
under the natural projection $P = \Pe_1 \times \Pe_3 \ra \Pe_3$. This
yields
$$
N_{C_0|P} = \co_{C_0}^{\oplus3}.
$$
Since $N_{C_0|Z} = \co_{C_0}(-2)$, the exact sequence
$0 \ra N_{C_0|Z} \ra N_{C_0|P} \ra N_{Z|P} \Big|C_0 \ra 0$ implies that there
are 2 possibilities for $N_{Z|P} \Big| C_0$, namely
$$
N_{Z|P} \Big| C_0 = \left\{ \begin{array}{cc}
\co_{C_0} \oplus \co_{C_0}(2)&\\
& \mbox{or}\\
\co_{C_0}(1)^{\oplus 2} &
\end{array}
\right.
\eqno(10)
$$
\smallskip

\noindent
On the other hand  the line bundle $\cl^{-1}$ on $Z$ of Proposition 5.3
may be considered as the normal bundle  of $Z$ in $Y$. The exact
sequence of normal bundles for $Z \subseteq Y \subseteq P$ shows that
$\cl^{-1}$ is a subbundle of $N_{Z|P}$.\\
But
$$
\cl^{-1}|C_0 = \co_Z (-2 C_0)| C_0 = \co_{C_0} (4)
$$
which cannot be a subbundle of $N_{Z|P} \Big| C_0$ according to (10).
\hfill $\Box$
\vspace{0.61cm}
\newline
 So for the rest of the paper we may assume that $Y$ is a double
 structure on $Z = \Pe_1 \times \Pe_1$, embedded in $P$ as in Lemma 5.2.
 For the computation of the cohomology of $\ce$ it will be important to
 know for which pairs $(a,b)$ the scheme $Y$ is $(a,b)$--normal in $P$.
\vspace{0.61cm}
\newline
Recall that a closed subscheme $V \subseteq \Pe_1 \times \Pe_3$ is
called $(a,b)$--{\it normal} for some nonnegative integers $a$ and $b$ if the
canonical map $H^0(\co_{\Pe_1 \times \Pe_3} (a,b)) \ra H^0(\co_V (a,b))$
is surjective.
 Note first\\
\bigskip

\noindent
{\bf Lemma 5.5} {\it Let $ i = (\vp_1, \vp_3) : Z \ra
\Pe_1 \times \Pe_3$ be the embedding of Lemma } 5.2. {\it Then }
$Z$ {\it is
$(a,b)$--normal in $\Pe_1 \times \Pe_3$ for all} $(a,b) \geq 0$.\\
\bigskip

\noindent
{\bf Proof.} Consider the commutative diagram
$$
\begin{array}{ccc}
H^0( \co_{\Pe_1 \times \Pe_3}(a,b)) &
\stackrel{(\vp_1, \vp_3)^{\ast}}{\makebox[1.2cm]{\rightarrowfill}}  & H^0(\co_Z
(a,b))\\[2ex]
(\vp_1 \times \vp_3)^{\ast} \searrow && \nearrow \Delta^{\ast}\\[2ex]
& H^0(Z \times Z, \vp^{\ast}_1 \co_{\Pe_1} (a) \otimes \vp_3^{\ast}
\co_{\Pe_3} (b)) &
\end{array}
$$
where $\Delta : Z \ra Z \times Z$ denotes the diagonal map. It suffices
to show  that for all $a, b \geq 0:$
\begin{itemize}
\item[(i)] $ ( \vp_1 \times \vp_3)^{\ast}$ is surjective
\item[(ii)]
$\Delta^{\ast}$ is surjective.
\end{itemize}
\underline{(i)}: By   K\"unneth $(\vp_1 \times \vp_3)^{\ast}$ identifies
with
$$
\vp^{\ast}_1 \otimes \vp^{\ast}_3 : H^0 (\co_{\Pe_1} (a)) \otimes H^0
(\co_{\Pe_3} (b)) \ra H^0 (\co_Z (af) \otimes H^0(\co_Z (b C_0 +bf).
$$
But $\vp^{\ast}_1$ is an isomorphism and $\vp^{\ast}_3$ is surjective
since $\vp_3$ is the Veronese embedding $\vp_3 : Z = \Pe_1 \times \Pe_1
\ra \Pe_3$,
which is projectively normal. \\
\underline{(ii)}: Let $\Pe^i_1$ denote the $i$--th component of $Z =
\Pe_1 \times \Pe_1$. According to K\"unneth the following diagram commutes
$$
\begin{array}{ccc}
H^0(\co_Z (af))\otimes H^0(\co_Z(b C_0 +bf) & \stackrel{\Delta^{\ast}}{
\longrightarrow} &
H^0 (\co_Z (b C_0 + (a+b)f))\\[1ex]
|| && || \\[1ex]
H^0(\co_{\Pe^1_1} (a)) \otimes H^0(\co_{\Pe^1_1} (b) \otimes
H^0(\co_{\Pe^2_1} (b)) & \longrightarrow & H^0(\co_{\Pe^1_1}
(a+b)) \otimes H^0 (\co_{\Pe^2_1} (b))
\end{array}
$$
But this is surjective, since the multiplication map
$H^0(\co_{\Pe_1} (a)) \otimes H^0(\co_{\Pe_1} (b)) \lra H^0
(\co_{\Pe_1} (a+b))$ is surjective for all $a, b \geq 0$. \hfill
$\Box$\\
\bigskip

\noindent
{\bf Proposition 5.6.}
{\it In the case $e=0$ the embedding $Y \ra \Pe_1 \times \Pe_3$ is
$(a,b)$--normal for all} $a \geq 0, b \geq 3$ or $(a,b) = (0,1).$\\
\bigskip

\noindent
{\bf Proof.} Using $(4)$ we have the following commutative diagram with
exact rows
$$
\begin{array}{lccccccccr}
0 \ra & I_{Z| \Pe_1 \times \Pe_3}(a,b) & \ra & \co_{\Pe_1 \times \Pe_3}
(a,b) & \ra & \co_Z(a,b) & \ra 0\\[1ex]
& \downarrow && \downarrow && || &\\[1ex]
0 \ra & \cl(a,b) & \ra & \co_Y (a,b) & \ra & \co_Z(a,b) & \ra 0.
\end{array}
$$
Since $Z$ is $(a,b)$--normal in $\Pe_1 \times  \Pe_3$, this yields
$$
\begin{array}{lcccccr}
0 \ra & H^0(I_{Z|\Pe_1 \times \Pe_3} (a,b)) & \ra & H^0(\co_{\Pe_1 \times
\Pe_3} (a,b)) & \ra & H^0 (\co_Z (a,b)) & \ra 0 \\[1ex]
&\downarrow \vp & & \downarrow && ||\\[1ex]
0 \ra & H^0(\cl(a,b)) & \ra & H^0(\co_Y (a,b)) & \ra & H^0(\co_Z(a,b)) &
\ra 0.
\end{array}
$$
Hence it suffices to show that $\vp$ is surjective. But
$$
\cl(a,b) = \co_Z ((b+2) C_0 + (a+b -2) f).
$$
So $\vp$ is surjective if and only if the canonical map
$$
H^0(N^{\ast}_{Z| \Pe_1 \times \Pe_3} (b C_0 +( a+b)f)) \ra H^0(\co_Z
((b+2) C_0 +(a+b -2)f)
$$
is surjective, where $N^{\ast}_{Z| \Pe_1 \times \Pe_3}$ denotes the
conormal bundle of $Z$ in $\Pe_1 \times \Pe_3$. Now
det$N^{\ast}_{Z| \Pe_1 \times \Pe_3} = \co_Z(-2 C_0 - 4f)$ and
the normal bundle  sequence for $Z \subset Y \subset P$ gives
$$
\scriptstyle{0 \ra \co_Z((b-4) C_0 +(a+b-2)f) \ra N^{\ast}_{Z| \Pe_1 \times
\Pe_3} (b C_0 +(a+b)f)) \ra \co_Z((b+2)C_0 + (a+b-2)f) \ra 0 .}
$$
But
$$
\begin{array}{rcl}
h^1(\co_Z (b-4) C_0 +(a+b - 2)f)&  = &h^1(\co_{\Pe_1}(b-4)) \cdot
h^0(\co_{\Pe_1}(a+b-2))\\
&+&h^0 (\co_{\Pe_1} (b-4)) \cdot h^1(\co_{\Pe_1} (a+b-2))\\
&=& 0 \quad \mbox{if} \quad b \geq 3 \quad \mbox{or} \quad  (a,b) = (0,1).
\hspace{2cm} \Box
\end{array}
$$
\bigskip

\noindent
{\bf Remark 5.7 } The proof of Proposition 5.6 yields
\\
for $b=2, a \geq 0: \quad \codim (H^0(\co_{\Pe_1 \times \Pe_3} (a,2))
\ra H^0 (\co_Y (a,2)) \leq a+1$\\
for $ b=1, a \geq 0: \quad  \codim (H^0(\co_{\Pe_1 \times \Pe_3} (a,1))
\ra H^0(\co_Y (a,1)) \leq 2a \qquad $ and\\
for $b =0, a \geq 1: \quad \codim (H^0(\co_{\Pe_1 \times \Pe_3} (a,0))
\ra H^0(\co_Y (a,0) \leq 3 (a-1)$.
\section{Splitting type of $\ce$} There are 2 types of lines in $\Pe_1
\times \Pe_3$. In this section lines $\{ t \} \times \ell $ (respectively
$\Pe_1
\times \{ x \}$) are called {\it horizontal} (respectively {\it
vertical}) {\it lines} in $\Pe_1 \times \Pe_3$. Restriction of the
vector bundle $\ce$ to these lines leads to the notion of {\it
horizontal} (respectively {\it vertical}) {\it generic splitting types}
as well as {\it
horizontal} (respectively {\it vertical}) {\it jumping lines.}
\bigskip

\noindent
{\bf Proposition 6.1 } (a) {\it The horizontal generic splitting type of
$\ce$ is } $(2,2)$.\\
(b) {\it The vertical generic splitting type of $\ce$ is} $(1,1)$.
\bigskip

\noindent
{\bf Proof.} (a) follows from the fact that $\ce_t$ is semistable but not
stable with $c_1 (\ce_t) =4$ on $\Pe_3$ (see Section 4).\\
(b): According to Section 1(4) the projection $X \subseteq \Pe_1 \times
\Pe_3 \ra \Pe_3$ is birational onto its image. Hence there is a a vertical
line  $\ell = \Pe_1 \times \{ x \}$ intersecting $X$ transversally in
exactly one point. Restricting the exact sequence Section 2(1) to $\ell$
gives
$$
0 \ra \co_{\ell} \ra \ce| \ell \ra \co_{\ell}/I_{\ell \cap A} \otimes
\co_{\ell} (2) \ra 0.
$$
Hence $\ce | \ell = \co_{\ell} (1) \otimes \co_{\ell} (1)$. \hfill $\Box$
\bigskip

\noindent
{\bf Proposition 6.2.} (a) {\it The horizontal jumping lines of $\ce$
are exactly the horizontal lines $\ell$ in $\Pe_1 \times \Pe_3$
intersecting the surface $Z \subset \Pe_1 \times \Pe_3$ of Proposition}
5.1.\\
(b) {\it If $\ell$ intersects $Z$ transversally, then} $\ce|\ell =
\co_{\ell} (4) \oplus \co_{\ell} (0)$.\\

\par
In particular all horizontal jumping lines are "higher" jumping
lines.
\bigskip

\noindent
{\bf Proof.} Let $\ell = \{ t \} \times \Pe_1$ intersect the line $Z_t$
in $\{ t \} \times \Pe_3$ transversally. Then $\ell$ intersects the
double structure $Y_t$ on $Z_t$ with multiplicity 2. Restricting the
exact sequence (2) of Section 4 to $\ell$ we get
$$
0 \ra \co_{\ell} \ra E_t (-2)| \ell \ra \co_{\ell}/I_{\ell \cap Y_t}
\oplus \co_{\ell} \ra 0
$$
which fills up to an exact sequence
$$
0 \ra \co_{\ell} (2) \ra E_t (-2) | \ell \ra \co_{\ell} (-2) \ra 0
$$
which necessarily splits. This proves (b).
\vspace{1ex}
\newline
{\it As for (a):} Since $\ce_t$ is semistable on $\Pe_3$,
according to a theorem of Barth (see [O.S.S] p. 228) the jumping lines
of $\ce_t$ on $\Pe_3$ form a divisor in the Grassmannian Gr$(1,3)$ and
there is an effective divisor  $D_{\ce_t}$ of degree
$c_2(\ce_t(-2))=2$ with support the set of jumping lines $S_{\ce_t}$ in $Gr
(1,3)$.\\
On the other hand: The hypersurface $H_{Z_t} = \{ \ell \in  Gr(1,3) |
\ell \cap Z_t = \emptyset \}$ consists of jumping lines. Since every
such jumping line is "higher", we get $D_{E_t} = 2H_{Z_t}$. This implies
the assertion. \hfill $\Box$\\

\par
Finally let us determine the vertical jumping lines of $\ce$. In order to state
the result recall from [L] that the coordinates of $\Pe_3$ can be chosen in
such
a way that the map $\vp_3 : X \ra \vp_3 (X) \subseteq \Pe_3$ is an
embedding outside the coordinate planes, whereas $\vp_3(X)$ is singular
along the coordinate planes. Choosing the coordinates of $\Pe_3$ in this
way, we have
\bigskip

\noindent
{\bf Proposition 6.3.} {\it The vertical jumping lines are exactly the
lines $\Pe_1 \times \{ x \}$ with $x$ contained in a coordinate
plane.\\
(b) For a general vertical jumping line }$\ell : \ce | \ell = \co_{\ell}
(2) \oplus  \co_{\ell} (0)$.
\bigskip

\noindent
For the proof let $S_{\ce} \subseteq \Pe_3$ denote the set of jumping lines
and consider the projection $q : \Pe_1 \times \Pe_3 \ra \Pe_3$.
\bigskip

\noindent
{\bf Lemma 6.4.} $S_{\ce} = \mbox{supp} (R^1_{q_{\ast}} \ce (-2, -2)).$
\bigskip

\noindent
{\bf Proof.} Let $U \subset \Pe_3$ be open affine. For $x \in U$ denote
$\ell = \Pe_1 \times \{ x \}$ and $I_{\ell}$ the ideal sheaf of $\ell$.
Then $H^2 (\Pe_1 \times U, \ce (-2, -2) \otimes I_{\ell}) =0$, since
$\Pe_1 \times U$ can be covered by 2 open affine sets. Hence the
canonical map $H^1(\Pe_1 \times U, \ce(-2, -2) ) \ra H^1(\ell, \ce (-2,
-2)|\ell)$ is surjective. This implies that the base change homomorphism
$$
R^1 q_{\ast} \ce (-2, -2) (x) \ra H^1 (\ell, \ce(-2, -2)|\ell)
$$
is an isomorphism. On the other hand $\ell = \Pe_1 \times \{ x \}$ is a
vertical jumping line if and only if $h^1(\ell, \ce (-2, -2)| \ell) >
0$. \hfill $\Box$\\

\par
In order to define a scheme structure on $S_{\ce}$ choose a resolution
$$
0 \ra \ck \ra \bigoplus^{r+2}_{i=1} \co_{\Pe_1 \times \Pe_3} (-a_i, -b_i)
\ra \ce (-2, -2) \ra 0 \eqno(1)
$$
with $a_i, b_i > 0$ for all $i$ and $\ck$ is a locally free sheaf of
rank $r$ on $\Pe_1 \times \Pe_3$. The sequence
$$
0 \ra R^1_{q_{\ast}} \ck \stackrel{\vp}{\ra} \bigoplus^{r+2}_{i=1}
R^1_{q_{\ast}} \co(-a_i, -b_i) \stackrel{\psi}{\ra} R^1_{q_{\ast}}
\ce(-2, -2) \ra 0 \eqno(2)
$$
is exact: $\vp$ is injective since $q_{\ast} \ce(-2, -2)$ is torsion
free and 0 outside $S_{\ce}$ and $\psi$ is surjective since
 $R^2_{q_{\ast}} \ck =0$.
\par
The sheaves $R^1_{q_{\ast}} \ck$ and $R^1_{q_{\ast}} \co(-a_i, -b_i)$ are
locally free
by the base change theorem since $h^1( \co_{\ell} (-a_i, -b_i))$ and
$h^1(\ck| \ell) = -\chi(\ck| \ell)$ are independent of $\ell$.
\par
Hence $\vp$ is a homomorphism of locally free sheaves of the same rank on
$\Pe_3$ and
$$
J_{\vp} := \mbox{Im} (\det \vp) \otimes \det \left(
\bigoplus^{r+2}_{i=1} R^1_{q_{\ast}} \co (-a_i, -b_i)\right)^{\ast} \subseteq
\co_{\Pe_3}
$$
is an invertible ideal sheaf in $\co_{\Pe_3}$ with supp$(\co_{\Pe_3} /
J_{\vp}) = S_{\ce}$.
Hence $D_{\ce} := (S_{\ce}, \co_{\Pe_3}/J_{\vp})$ is a divisor on
$\Pe_3$.
\bigskip

\noindent
{\bf Lemma 6.5.} deg$D_{\ce} =4$.\\
\medskip

\noindent
{\bf Proof.}  Let  $\Pe_1 \subseteq \Pe_3$ be a general line. Since the
divisor $D_{\ce}$ intersects the line $\Pe_1$ in a divisor of the same
degree, we may restrict the whole situation to $\Pe_1 \times \Pe_1$ and
compute the degree of the corresponding divisor on $\Pe_1$.
\par
By abuse of notation we denote the restricted objects to $\Pe_1 \times
\Pe_1$ by the same letter. In particular we have the restricted
sequences (1) on $\Pe_1 \times \Pe_1$ and (2) on $\Pe_1$.\\
If $h_1 = [p^{\ast} \co_{\Pe_1} (1)]$ and $h_2 = [q^{\ast} \co_{\Pe_1}
(1)]$, then from (1) we get
$$
c_1 (K)= - \sum_i (a_i h_1 + b_i h_2) + 2 h_1 \eqno(3)
$$
$$
c_2 (K) = \sum_{i \not=j} a_i b_j - 4 -2 \sum_i b_i \eqno(4)
$$
Applying flat base change and the projection formula, we get
$$
\begin{array}{rl}
R^1 q_{\ast} \co (-a_i, -b_i) & = R^1_{q_{\ast} } p^{\ast} \co_{\Pe_1}
(- a_i) \otimes
\co_{\Pe_1} (- b_i)\\[1ex]
&= H^1 (\co_{\Pe_1} (-  a_i)) \otimes
\co_{\Pe_1} (-b_i)
\end{array}
$$
and hence
$$
c_i (R^1 q_{\ast} \co(-a_i, -b_i)) =-b_i (a_i - 1).
$$
The Theorem of Grothendieck--Riemann--Roch for the morphism $q: \Pe_1
\times
\Pe_1 \ra \Pe_1$ and equations (3) and (4) give
$$
\begin{array}{ll}
c_1(R^1 q_{\ast} \ck) &= - c_1 (\ck)h_1 - \frac{1}{2} (c^2_1 (\ck) -
2c_2(\ck))\\[1ex]
& = -
\D{ \sum_i a_i b_i + \sum_i b_i - 4}
\end{array}
$$
It follows that
$$
\begin{array}{ll}
\deg D_{\ce} & = - c_1 \left( \bigoplus^{r+2}_{i=1} R^1  q_{\ast}
\co(-a_i, -b_i)\right) - c_1(R^1_{q_{\ast}} \ck)\\[1ex]
& = \D{- \sum_i b_i (a_i -1) + \sum_i a_i b_i - \sum_i b_i +4 = 4}
\hspace{5.0cm}
\Box
\end{array}
$$
{\bf Proof of Proposition 6.3.} Let $H_1, \ldots , H_4$ denote the 4
coordinate planes in $\Pe_3$. For a general point $x$ of $\vp_3 (X) \cap
H_i$ the vertical line $\ell = \Pe_1 \times \{ x \}$ intersects
$X \subset \Pe_1 \times \Pe_3$ transversally in 2 points. Restriction of
the exact sequence
(1) of Section 2 to $\ell$ gives $\ce | \ell = \co_{\ell}(2) \oplus
\co_{\ell} (2)$. Hence the divisor $\vp_3(X) \cap H_i \subset D_{\ce}$
for $i = 1, \ldots , 4$. Since $H_1 + \ldots + H_4$ is the only divisor
of degree $\leq 4$ in $\Pe_3$ containing the 4 curves $\vp_3 (X) \cap
H_i$, the assertion follows. \hfill $\Box$
\section{Problems}
Of course concerning the vector bundle $\ce$ on $\Pe_1 \times \Pe_3$ one
can ask all the questions which have been studied for the
Horrocks--Mumford
bundle on $\Pe_4$. The most important problems are
\bigskip

\noindent
(1) {\bf Cohomology of $\ce$:} Compute the dimensions of the
groups $H^i(\ce(a,b))$ for all $a,b$ and $i$. I computed most of these
groups, unfortunately  the list is not complete. In particular the most
important dimension $h^0(\ce)$ is not known.
\bigskip

\noindent
(2) {\bf Automorphisms of $\ce$:} According to Proposition 3.5 the group
of automorphisms of the vector bundle $\ce$ contains the dihedral group
$D_8 = H_1(L_1, L_3)$. Does $\ce$ admit further automorphisms? Recall
that the automorphism group of the Horrocks--Mumford bundle is the
normalizer of the Heisenberg group  in $S\ell_5(\C)$. Is the analogue
statement valid in the case of $\ce$?
\bigskip

\noindent
(3) {\bf Moduli:} Describe the moduli spaces of vector bundles with the
invariants of $\ce$. This question is closely related to (1) and (2).
For example, if $h^0(\ce) \leq 2$, then $\ce$ admits moduli, since there
is a two--dimensional family of abelian surfaces in $\Pe_1 \times
\Pe_3$.
\bigskip

\noindent
(4) {\bf Other constructions:} It seems not too difficult to construct
the vector bundle $\ce$ out of the double structure $Y$ on $\Pe_1 \times
\Pe_1$ in $\Pe_3 \times \Pe_3$ of Section 5. Are there other
constructions?
\bigskip

\noindent
(5) {\bf Degenerations:} It is not difficult to see that every smooth
zero set of a global section of $\ce$ is an abelian surface $A$ in
$\Pe_1 \times \Pe_3$. Determine the degenerations of $A$ in $\Pe_1
\times \Pe_3$. Is every such degeneration the zero set of a section of
(one of the vector bundles) $\ce$?
\bigskip

\noindent
(6) {\bf Classification:} Are there other rank--2 vector bundles on $\Pe_1
\times \Pe_3$ apart
from the bundles derived from $\ce$ or the bundles of the introduction?
\vspace{2cm}
\newline
{\bf Literature:}
\begin{quote}\begin{itemize}
\item[\mbox{\makebox[1.2cm][l]{[CAV]}}] {\bf H. Lange, Ch. Birkenhake:} {\it
Complex Abelian Varieties};
Grundlehren 302, Springer 1992
\item[\mbox{\makebox[1.2cm][l]{[H]}}] {\bf R. Hartshorne:} {\it Stable Vector
Bundles of Rank 2 on} $\Pe_3$;
Math. Ann. 238 (1978), 229 -280
\item[\mbox{\makebox[1.2cm][l]{[L]}}] {\bf H. Lange:} {\it Abelian surfaces in
}$\Pe_1 \times \Pe_3$.
Arch. Math. 63 (1994), 80 -- 84.
\item[\mbox{\makebox[1.2cm][l]{[O.S.S]}}] {\bf Ch. Okonek, M. Schneider, H.
Spindler:} {\it
Vector Bundles on Complex  Projective Spaces}; Progr. in Math. 3,
Birkh\"auser 1980
\end{itemize}
\end{quote}
\end{document}